\documentclass[12pt,a4paper,twoside]{article}

\newcommand{\pageformat}[6]{\setlength{\hoffset}{-1in}
                  \setlength{\voffset}{-1in}
                  \addtolength{\hoffset}{#5}
                            \addtolength{\voffset}{#6}
                            \setlength{\oddsidemargin}{#1}
                            \setlength{\evensidemargin}{#2}
                            \setlength{\textwidth}{\paperwidth}
                  \addtolength{\textwidth}{-\oddsidemargin}
                  \addtolength{\textwidth}{-\evensidemargin}
                  \addtolength{\textwidth}{-\marginparsep}
                  \addtolength{\textwidth}{-\marginparwidth}
                            \setlength{\topmargin}{#3}
                            \setlength{\textheight}{\paperheight}
                  \addtolength{\textheight}{-\topmargin}
                  \addtolength{\textheight}{-\headheight}
                  \addtolength{\textheight}{-\headsep}
                  \addtolength{\textheight}{-\footskip}
                  \addtolength{\textheight}{-#4}}
\pageformat{2cm}{3cm}{25mm}{25mm}{1pt}{0pt}

\usepackage{ifthen}
\newboolean{article}
    \setboolean{article}{true}
\newboolean{report}
\newboolean{book}
\newboolean{letter}
\newboolean{german}
\newboolean{italian}
\newboolean{nobaselinestretch}
\newboolean{nosectionappendix}
\newboolean{oldtoc}
\newboolean{nosectionequn}
\newboolean{notheorem}
\ifthenelse{\boolean{german}}{
    \usepackage{german}}{}

\usepackage[latin1]{inputenc}
\usepackage{amsmath}
\usepackage{amssymb}
\usepackage[mathscr]{eucal}

\ifthenelse{\boolean{notheorem}}{}{
    \usepackage{theorem}}

\ifthenelse{\boolean{nobaselinestretch}}{}{
    \renewcommand{\baselinestretch}{1.25}}

\newenvironment{env}[2]{\begin{#1}#2\end{#1}}{}
    \newcommand{\beq}[1]{\begin{env}{equation}{#1}}
    \newcommand{\beqn}[1]{\begin{env}{equation*}{#1}}
    \newcommand{\bal}[1]{\begin{env}{align}{#1}}
    \newcommand{\baln}[1]{\begin{env}{align*}{#1}}
    \newcommand{\bga}[1]{\begin{env}{gather}{#1}}
    \newcommand{\bgan}[1]{\begin{env}{gather*}{#1}}
    \newcommand{\bflal}[1]{\begin{env}{flalign}{#1}}
    \newcommand{\bflaln}[1]{\begin{env}{flalign*}{#1}}
    \newcommand{\bmu}[1]{\begin{env}{multline}{#1}}
    \newcommand{\bmun}[1]{\begin{env}{multline*}{#1}}
    \newcommand{\bsp}[1]{\begin{env}{split}{#1}}

    \newcommand{\eeq}{\end{env}}
    \newcommand{\eeqn}{\end{env}}
    \newcommand{\eal}{\end{env}}
    \newcommand{\ealn}{\end{env}}
    \newcommand{\ega}{\end{env}}
    \newcommand{\egan}{\end{env}}
    \newcommand{\eflal}{\end{env}}
    \newcommand{\eflaln}{\end{env}}
    \newcommand{\emu}{\end{env}}
    \newcommand{\emun}{\end{env}}
    \newcommand{\esp}{\end{env}}

\newcommand{\lf}{\vspace{2ex}}

\renewcommand{\bf}[1]{\textbf{#1}}
\renewcommand{\it}[1]{\textit{#1}}

\renewcommand{\sf}[1]{\textsf{#1}}

\renewcommand{\tt}[1]{\texttt{#1}}
\newcommand{\hl}[1]{\bf{\it{#1}}}

\newcommand{\msf}[1]{\text{\small $\sf{#1}$}}

\newcommand{\cmc}[1]{\mathcal{#1}}
\newcommand{\eus}[1]{\mathscr{#1}}
\newcommand{\euf}[1]{\mathfrak{#1}}
\newcommand{\bb}[1]{\mathbb{#1}}
\newcommand{\msmall}[1]{{\setlength{\arraycolsep}{.6ex}\text{\small$#1$}}}

\newcommand{\mscriptsize}[1]{{\setlength{\arraycolsep}{.3ex}\text{\scriptsize$#1$}}}
\newcommand{\mtiny}[1]{{\setlength{\arraycolsep}{.3ex}\text{\tiny$#1$}}}
\newcommand{\nbd}[1]{$#1$\nobreakdash--}
\newcommand{\ol}[1]{\overline{#1}}
\newcommand{\ul}[1]{\underline{#1}}

\newcommand{\abs}[1]{\left\lvert#1\right\rvert}
\newcommand{\norm}[1]{\left\lVert#1\right\rVert}

\newcommand{\bnorm}[1]{\bigl\lVert#1\bigr\rVert}

\newcommand{\Bnorm}[1]{\Bigl\lVert#1\Bigr\rVert}

\newcommand{\bfam}[1]{\bigl(#1\bigr)}
\newcommand{\Bfam}[1]{\Bigl(#1\Bigr)}
\newcommand{\AB}[1]{\langle#1\rangle}

\newcommand{\BAB}[1]{\Bigl\langle#1\Bigr\rangle}
\newcommand{\CB}[1]{\{#1\}}
\newcommand{\bCB}[1]{\bigl\{#1\bigr\}}
\newcommand{\BCB}[1]{\Bigl\{#1\Bigr\}}
\newcommand{\SB}[1]{[#1]}
\newcommand{\bSB}[1]{\bigl[#1\bigr]}

\newcommand{\Matrix}[1]{\begin{pmatrix}#1\end{pmatrix}}
\newcommand{\SMatrix}[1]{\msmall{\Matrix{#1}}}

\newcommand{\tMatrix}[1]{\mtiny{\Matrix{#1}}}
\newcommand{\rtMatrix}[1]{\raisebox{.3ex}{\tMatrix{#1}}}

\newcommand{\sbars}[1]{\:\bar{#1}^s\:}
\newcommand{\sodots}{\sbars{\odot}}
\newcommand{\set}[2][]{
    \ifthenelse{\equal{#1}{}}{
        \CB{#2}}{
        \CB{#1~|~#2}}}
\newcommand{\bset}[2][]{
    \ifthenelse{\equal{#1}{}}{
        \bCB{#2}}{
        \bCB{#1~|~#2}}}
\newcommand{\Bset}[2][]{
    \ifthenelse{\equal{#1}{}}{
        \BCB{#2}}{
        \BCB{#1~\big|~#2}}}
\newcommand{\zero}{\CB{0}}

\DeclareMathOperator{\ls}{\normalfont\msf{span}}

\DeclareMathOperator{\fls}{\mscriptsize{\sf{span}}}
\DeclareMathOperator{\cls}{\ol{\ls}}

\DeclareMathOperator{\id}{\normalfont\msf{id}}

\renewcommand{\ker}{\operatorname{\msf{ker}}}

\newcommand{\C}{\bb{C}}

\newcommand{\E}{\bb{E}}

\newcommand{\cA}{\cmc{A}}
\newcommand{\cB}{\cmc{B}}
\newcommand{\cC}{\cmc{C}}

\newcommand{\sB}{\eus{B}}

\newcommand{\sD}{\eus{D}}

\newcommand{\sL}{\eus{L}}

\newcommand{\sN}{\eus{N}}

\newcommand{\eD}{\euf{D}}

\newcommand{\eH}{\euf{H}}

\newcommand{\eK}{\euf{K}}
\newcommand{\eL}{\euf{L}}

\ifthenelse{\boolean{nosectionequn}}{}{
    \numberwithin{equation}{section}
    }

\ifthenelse{\boolean{article}\or\boolean{letter}\or\boolean{nosectionequn}}{
    \setboolean{nosectionappendix}{true}}{}
\ifthenelse{\boolean{nosectionappendix}}{}{
    \renewcommand{\appendix}{
        \chapter*{\appendixname}
        \addcontentsline{toc}{chapter}{\appendixname}
        \renewcommand{\thesection}{\Alph{section}}
        \setcounter{section}{0}}}
   
\ifthenelse{\boolean{report}\or\boolean{book}}{
    }{}

\ifthenelse{\boolean{notheorem}}{}{
        \newcommand{\mnname}{Mathematical note.}
        \newcommand{\enname}{End of the note.}
        \newcommand{\definame}{Definition.}
        \newcommand{\propname}{Proposition.}
        \newcommand{\lemname}{Lemma.}
        \newcommand{\exname}{Example.}
        \newcommand{\exername}{Exercise.}
        \newcommand{\remname}{Remark.}
        \newcommand{\obname}{Observation.}
        \newcommand{\thmname}{Theorem.}
        \newcommand{\corname}{Corollary.}
        \newcommand{\proofname}{Proof.}
    \ifthenelse{\boolean{german}}{
        \renewcommand{\mnname}{Mathematische Notiz.}
        \renewcommand{\enname}{Ende der Notiz.}
        \renewcommand{\exname}{Beispiel.}
        \renewcommand{\exername}{Übung.}
        \renewcommand{\remname}{Bemerkung.}
        \renewcommand{\obname}{Beobachtung.}
        \renewcommand{\thmname}{Satz.}
        \renewcommand{\corname}{Korollar.}
        \renewcommand{\proofname}{Beweis.}}{}
    \ifthenelse{\boolean{italian}}{
        \renewcommand{\mnname}{Nota matematica.}
        \renewcommand{\enname}{Fina della nota.}
        \renewcommand{\definame}{Definizione.}
        \renewcommand{\propname}{Proposizione.}
        \renewcommand{\exname}{Esempio.}
        \renewcommand{\exername}{Esercizio.}
        \renewcommand{\remname}{Nota.}
        \renewcommand{\obname}{Osservazione.}
        \renewcommand{\thmname}{Teorema.}
        \renewcommand{\corname}{Corollario.}
        \renewcommand{\proofname}{Dimostrazione.}

       \renewcommand{\appendixname}{Appendice}

       }{}
    \theoremheaderfont{\normalfont\bfseries}
    \theoremstyle{change}
        \theorembodyfont{\rmfamily}
            \newtheorem{emp}{}[section]
                \newcommand{\bemp}[1][]{
                    \begin{emp}\hskip-\labelsep\bf{#1}\hskip\labelsep}
                \newcommand{\eemp}{\end{emp}}
\newtheorem{itemp}[emp]{}
                \newcommand{\bitemp}[1][]{
                    \begin{itemp}\hskip-\labelsep\bf{#1}\hskip\labelsep\normalfont\itshape}
                \newcommand{\eitemp}{\end{itemp}}
            \newtheorem{mn}[emp]{\mnname}
                \newcommand{\bnm}{\begin{mn}~\begin{quotation}\renewcommand{\baselinestretch}{1}\small\noindent\ignorespaces}
                \newcommand{\enm}{\end{quotation}\hfill\bf{\enname}\end{mn}}
            \newtheorem{ex}[emp]{\exname}
                \newcommand{\bex}{\begin{ex}}
                \newcommand{\eex}{\end{ex}}
            \newtheorem{exer}[emp]{\exername}
                \newcommand{\bexer}{\begin{exer}}
                \newcommand{\eexer}{\end{exer}}
            \newtheorem{defi}[emp]{\definame}
                \newcommand{\bdefi}{\begin{defi}}
                \newcommand{\edefi}{\end{defi}}
            \newtheorem{rem}[emp]{\remname}
                \newcommand{\brem}{\begin{rem}}
                \newcommand{\erem}{\end{rem}}
            \newtheorem{ob}[emp]{\obname}
                \newcommand{\bob}{\begin{ob}}
                \newcommand{\eob}{\end{ob}}
        \theorembodyfont{\normalfont\itshape}
            \newtheorem{thm}[emp]{\thmname}
                \newcommand{\bthm}{\begin{thm}}
                \newcommand{\ethm}{\end{thm}}
            \newtheorem{prop}[emp]{\propname}
                \newcommand{\bprop}{\begin{prop}}
                \newcommand{\eprop}{\end{prop}}
            \newtheorem{cor}[emp]{\corname}
                \newcommand{\bcor}{\begin{cor}}
                \newcommand{\ecor}{\end{cor}}
            \newtheorem{lem}[emp]{\lemname}
                \newcommand{\blem}{\begin{lem}}
                \newcommand{\elem}{\end{lem}}
\newenvironment{empn}[1]{\lf\noindent\bf{#1}\ignorespaces\hskip\labelsep}{\lf}
		\newcommand{\bempn}[1]{\begin{empn}{#1}}
		\newcommand{\eempn}{\end{empn}}
		\newcommand{\bitempn}[1]{\begin{empn}{#1}\normalfont\itshape}
		\newcommand{\eitempn}{\end{empn}}
                \newcommand{\bnmn}{\begin{empn}{\mnname}~\begin{quotation}\renewcommand{\baselinestretch}{1}\small\noindent\ignorespaces}
                \newcommand{\enmn}{\end{quotation}\hfill\bf{\enname}\end{empn}}
		\newcommand{\bexn}{\begin{empn}{\exname}}
		\newcommand{\eexn}{\end{empn}}
		\newcommand{\bexern}{\begin{empn}{\exername}}
		\newcommand{\eexern}{\end{empn}}
		\newcommand{\bdefin}{\begin{empn}{\definame}}
		\newcommand{\edefin}{\end{empn}}
		\newcommand{\bremn}{\begin{empn}{\remname}}
		\newcommand{\eremn}{\end{empn}}
		\newcommand{\bobn}{\begin{empn}{\obname}}
		\newcommand{\eobn}{\end{empn}}

\newcommand{\qedsymbol}{~\rule[-0.35mm]{2mm}{2mm}}
    \newcounter{proof}[emp]
    \newenvironment{Proof}[1]{
        \vspace{1ex}
        \renewcommand{\item}[1][\stepcounter{proof}(\roman{proof})]%
            {##1\hskip\labelsep}
        \noindent\textsc{#1\hskip\labelsep}}{
        \nolinebreak\qedsymbol}
    \newcommand{\proof}[1][\proofname]{
        \begin{Proof}{#1}\ignorespaces}
    \newcommand{\qed}{\end{Proof}}
    \newcommand{\noqed}{
        \renewcommand{\qedsymbol}{}
        \end{Proof}}}
    \ifthenelse{\boolean{italian}}{
        \renewcommand{\proofname}{Dimostrazione.}}{}

\usepackage[varg]{txfonts}

\renewcommand{\thefootnote}{[\arabic{footnote}]}

\usepackage{hyperref}

\begin{document}

\bibliographystyle{amsalpha}

\title{Hilbert-von Neumann Modules\\\it{versus}\\Concrete von Neumann Modules}

\author{Michael Skeide{\renewcommand{\thefootnote}{}}, May 2012, revised July 2021}

\date{}

\maketitle

\vspace{-3ex}
\begin{abstract}
\noindent
Apart from presenting some new insights and results, one of our main purposes is to put some records in the development of von Neumann modules straight.

The von Neumann or \nbd{W^*}objects among the Hilbert (\nbd{C^*})modules are around since the first papers by Paschke (1973) \cite{Pas73} and Rieffel (1974) \cite{Rie74,Rie74a} that lift Kaplansky's setting (1953) \cite{Kap53} to modules over noncommutative \nbd{C^*}algebras. While the formal definition of \it{\nbd{W^*}modules} is due to Baillet, Denizeau, and Havet (1988) \cite{BDH88}, the one of \it{von Neumann modules} as strongly closed operator spaces started with Skeide (2000) \cite{Ske00b}. It has been paired with the definition of \it{concrete von Neumann modules} in Skeide (2006) \cite{Ske06b}. It is well-known that (pre-)Hilbert modules can be viewed as \it{ternary rings of operators} and in their definition of \it{Hilbert-von Neumann modules}, Bikram, Mukherjee, Srinivasan, and Sunder (2012) \cite{BMSS12} take that point of view.

We illustrate that a (suitably nondegenerate) \it{Hilbert-von Neumann module} \cite{BMSS12} is the same thing as a \it{strongly full concrete von Neumann module} \cite{Ske06b}. We also use this occasion to put some things in the papers \cite{Ske00b,Ske06b,BMSS12} right. We show that the tensor product of (concrete, or not) von Neumann correspondences \bf{is}, indeed, (a generalization of) the tensor product of Connes correspondences (claimed in Skeide (2008) \cite{Ske08a}), viewed in a way quite different from \cite{BMSS12}. We also present some new arguments that are useful even  for (pre-)Hilbert modules.
\end{abstract}

\section{Von Neumann modules: Comparison} \label{compSEC}

Let us right start with comparing the notion of \it{Hilbert-von Neumann modules} in the sense of Bikram, Mukherjee, Srinivasan, and Sunder (2012) \cite{BMSS12} with that of \it{concrete von Neumann modules} from Skeide (2006) \cite{Ske06b}. On the relation with the, still earlier, notion of \it{von Neumann modules} in Skeide (2000) \cite{Ske00b}, we comment in Remark \ref{crem}, especially, in Footnotes \ref{embFN} and \ref{vNcvNFN}.

Let $G$ and $H$ be Hilbert spaces. For $E\subset\sB(G,H)$ denote by $\SB{E}:=\cls^sE$ the strongly closed subspace of $\sB(G,H)$ generated by $E$.

Let $\cB$ be a von Neumann algebra acting nondegenerately on the Hilbert space $G$. Skeide \cite[Definition 2]{Ske06b} says, $E$ is a \hl{concrete von Neumann \nbd{\cB}module} if $E$ is strongly closed, if it is a (right) \nbd{\cB}submodule%
\footnote{
In \cite{Ske06b} we omitted to repeat that a module is closed under addition.
}
of $\sB(G,H)$ (that is, $E+E\subset E$ and $E\cB\subset E$), if $E^*E\subset\cB$, and if $\cls EG=H$. We say $E$ is \hl{strongly full} if $\cls^sE^*E=\cB$.

Bikram, Mukherjee, Srinivasan, and Sunder \cite{BMSS12} say on Page 50, $E$ is a \hl{von Neumann corner} if $E=\SB{E}\supset EE^*E$. It is \hl{nondegenerate} if $\cls EG=H$, $\cls E^*H=G$. They say in \cite[Definition 1.2(1)]{BMSS12}, a \hl{Hilbert-von Neumann module} over a von Neumann algebra $\cA$ is a von Neumann corner $E$ with a normal isomorphism $\pi$ from $\cA$ onto $\SB{E^*E}$.

Modulo nondegeneracy conditions, these definitions capture the same class of subspaces of $\sB(G,H)$:
\bthm \label{=thm}
Let $G$ and $H$ be Hilbert spaces. For the subset $E$ of $\sB(G,H)$ denote by $\cB$ the strong closure in $\sB(G)$ of the algebra generated by $E^*E$. Then the following are equivalent:
\begin{enumerate}
\item
$E$ is a nondegenerate von Neumann corner (so that, in particular, $\cB=\SB{E^*E}$).\vspace{-1ex}

\item
$E$ is a Hilbert-von Neumann module over $\cB$ for $\pi=\id_\cB$ satisfying $\cls EG=H$ and $\id_G\in\cB$.\vspace{-1ex}

\item
$E$ is a(n, automatically, strongly full) concrete von Neumann \nbd{\cB}module.
\end{enumerate}
Moreover, for an isomorphism $\pi$ from the von Neumann algebra $\cA$ onto the von Neumann algebra $\id_G\in\pi(\cA)\subset\sB(G)$,  the subset $E\subset\sB(G,H)$ is a Hilbert-von Neumann module over $\cA$ satisfying $\cls EG=H$ if and only if $E$ is a strongly full concrete von Neumann \nbd{\pi(\cA)}module. 
% \footnote{%In \cite{BMSS12}, it is not specified if $\pi$ is an isomorphism of \nbd{W^*}algebras or an isomorphism of von Neumann algebras.
% We understand $\pi$ as an isomorphism of von Neumann algebras. This includes, in particular, that the von Neumann algebra $\pi(\cA)=\cB$, respectively, are von Neumann algebras and, therefore, act nondegenerately on the Hilbert space $H$.}
\ethm

\proof
This is immediate from the observation that ~$\ls EE^*E=\ls\Bfam{E\ls(E^*E)}$.\qed

\brem \label{crem}~

\begin{enumerate}
\item
The statement in the proof is essentially the necessary observation to recover the well-known fact that Hilbert modules and \it{ternary rings of operators } (\it{TRO}) are the same thing; see the book \cite{BlLMe04} by Blecher and Le Merdy (in particular, \cite[Section 8.7]{BlLMe04}) and the references therein. (The references go back at least to the to the 1990s.) A \hl{TRO} is a (usually, norm closed) linear subspace $E$ of a \nbd{C^*}algebra (for instance, of $\sB(G,H)\subset\sB(G\oplus H)$) satisfying $EE^*E\subset E$. Of course, nothing changes under strong closure.

\item
Von Neumann modules have been defined in Skeide \cite{Ske00b}. A pre-Hilbert module $E$ over a von Neumann algebra $\cB$ is a \hl{von Neumann module} if its \it{extended linking algebra} $\rtMatrix{\cB&E^*\\E&\sB^a(E)}=\sB^a\rtMatrix{\cB\\E}$ (denoting by $\sB^a(E)$ the bounded adjointable operators on $E$) is a von Neumann algebra.%
\footnote{ \label{embFN}
It is to be noted that a von Neumann algebra $\cB$ is, as always(\bf{!}), acting nondegenerately as $\cB\subset\sB(G)$ on a Hilbert space $G$. Having said this, we recall the well-known (at least, since Rieffel \cite{Rie74a}) fact that every pre-Hilbert module over a concrete operator pre-\nbd{C^*}algebra $\sB\subset\sB(G)$ embeds canonically (and uniquely, up to obvious unitary equivalence) into $\sB(G,H)$, where $H:=E\odot G$ is the Hausdorff completion of the algebraic tensor product $E\,\ul{\otimes}\,G$ with the semiinner product determined by $\AB{x\otimes g,x'\otimes g'}:=\AB{g,\AB{x,x'}g'}$. Now $\sB^a\rtMatrix{\cB\\E}$ sits naturally in $\sB\rtMatrix{G\\H}$, and it is clear what it means for $\sB^a\rtMatrix{\cB\\E}$ to be a von Neumann algebra in $\sB\rtMatrix{G\\H}$.\vspace{1ex}
}
Directly after the definition, in \cite[Proposition 4.5]{Ske00b} it is stated that this property is equivalent to $E$ being strongly closed, that is, $E$ is a concrete von Neumann module.
\end{enumerate}
The definition of concrete von Neumann \nbd{\cB}module frees the definition of von Neumann \nbd{\cB}module from the obligation to first construct the (anyway, uniquely determined by $\cB\subset\sB(G)$ and $E$) embedding $E\subset\sB(G,H)$; see Footnote \ref{embFN}. If $E$ is given as $E\subset\sB(G,H)$, then Theorem \ref{=thm} tells that the definition in \cite{BMSS12} does not present a mentionable difference to the definition in \cite{Ske06b}.%
\footnote{ \label{vNcvNFN}
Calling this an ``apparently much simpler manner'' avoiding a ``two stages completion'' in the abstract of \cite{BMSS12} seems out of place. The ``two stages completion'' they are referring to is that of norm closure (to get a \nbd{C^*}module) and then strong closure. It is not necessary to go through \nbd{C^*}modules in the construction, and while \cite[Definition 4.4]{Ske00b} is really (for convenience, in order not to frighten people who are afraid of pre-Hilbert modules) formulated only for Hilbert modules, \cite[Definition 3.1.1]{Ske01} (which for some reason is the only work quoted in \cite{BMSS12}) is formulated for pre-Hilbert modules. As for the (Hausdorff) completion to obtain $H$,  in \cite[Sections 2+3]{BMSS12} quite similar constructions are performed. It would appear strange to ``allow'' such constructions for Stinespring type theorems and the tensor product of correspondences, but not to obtain the tensor product $E\odot G$. (Anyway, it is well-known (at least, from or Skeide \cite{Ske00b}, but probably earlier; see also Murphy \cite{Mur97}) how to obtain the Stinespring construction of a CP-map into $\cB\subset\sB(G)$ from the tensor product of Pascke's GNS-correspondence \cite{Pas73} $E$ and the representation space $G$ of $\cB$. In this context, while it is not clear to us why \cite[Section 2]{BMSS12} is restricted to \it{standard} Hilbert-von Neumann modules, a comparison with \cite[Section 7]{Ske00b} seems interesting.)\vspace{1ex}
}
\erem

\section{Self-duality}

Von Neumann or \nbd{W^*}modules can be characterized as those objects in the category of (pre-)Hilbert modules over a von Neumann algebra that are \it{self-dual}; a (pre-)Hilbert \nbd{\cB}module $E$ is \hl{self-dual} if every bounded right linear map $\Phi\colon E\rightarrow\cB$ has the form $y^*\colon x\mapsto\Phi(x)=\AB{y,x}$ for a suitable (obviously, unique) $y\in E$. Baillet, Denizeau, and Havet \cite{BDH88} take the property to be self-dual to define \nbd{W^*}modules.%
\footnote{
Attention, a few authors by \it{\nbd{W^*}module} mean just a Hilbert module over a \nbd{W^*}algebra. But this seems not quite enough to merit to be considered an object in a \nbd{W^*}category.\vspace{1ex}
}
For other definitions, of course, self-duality has to be a theorem.%
\footnote{
One should note that the various definitions may be compared not only as to what extent they allow to check easily if a given pre-Hilbert module $E$ over a von Neumann algebra is a von Neumann module or not; checking strong closure seems to have considerable advantages over checking self-duality. But they may also be compared as to what extent it is easy to obtain a von Neumann module from $E$ in case it is not yet. Paschke \cite{Pas73} has constructed a self-dual completion. As a space it is simply $\sB^r(E,\cB)$, the space of all bounded right linear $\Phi\colon E\rightarrow\cB$. But while it is actually not \it{a priori} easy to find that space, it is even much harder to compute the inner products of its elements. The definitions in \cite{Ske00b,Ske01,Ske06b,BMSS12} suggest: Take the strong closure of $E$ -- and it is obvious what the inner product of elements $x,y$ in this closure is, namely, $\AB{x,y}=x^*y$, and that it takes values in the von Neumann algebra $\cB$.\vspace{1ex}
}
Actually, to show that a definition of von Neumann module is equivalent to another one, may be done by showing that it is equivalent to self-duality.%
\footnote{
In \cite[Proposition 1.9]{BMSS12} it is shown that Hilbert-von Neumann modules are self-dual, and it is claimed earlier that this shows that this is equivalent to being a von Neumann module in the sense of \cite{Ske00b}. But, actually the backwards direction of this statement is missing; it would not be easy to obtain it without the procedure of embedding a (pre-)Hilbert module into $\sB(G,H)$ (see Footnote \ref{embFN}), which the authors of \cite{BMSS12} are so keen to avoid.

Anyway, it appears a remarkable idea to show equivalence of Hilbert-von Neumann modules and concrete von Neumann modules appealing to the (not so easy) self-duality, while we have the (fairly obvious) Theorem \ref{=thm} (and the easy observation that the procedure explained in Footnote \ref{embFN} transforms a von Neumann module into a concrete von Neumann module and the observation that a concrete von Neumann module is, obviously, a von Neumann module).
}

Fortunately, we have:

\bitemp[Theorem ({\protect\cite[Theorem 4.16]{Ske00b}}).] \label{sdthm}
Any von Neumann module is self-dual.
\eitemp

We also have (more or less \cite[Theorem 5.5]{Ske00b}) the following (more or less obvious, when taking into account Footnote \ref{embFN}) consequence:

\bcor \label{vNiffcor}
A (pre-)Hilbert module over a von Neumann algebra is a von Neumann module (or isomorphic to a concrete von Neumann module or isomorphic to a Hilbert-von Neumann module) if and only if it is self-dual.
\ecor

The proof of Theorem \ref{sdthm} presented in \cite{Ske00b} relies on three basic pieces of knowledge:
\begin{itemize}
\item
Polar decomposition (\cite[Proposition 2.10]{Ske00b}%
\footnote{
See \cite[Lemma 1.6]{BMSS12}.
}
): Every $x$ in a von Neumann \nbd{\cB}module $E$ can be written (uniquely) as $v\abs{x}$ where $v$ is a partial isometry in $E$ (satisfying $\ker\AB{v,v}=\ker\AB{x,x}\subset\cB\subset\sB(G)$).

\item
Existence of complete quasi-orthonormal systems (\cite[Theorem 4.11]{Ske00b}%
\footnote{
See also the proof of \cite[Proposition 1.7]{BMSS12}.
}
): In a von Neumann \nbd{\cB}module $E$ there exists a family $\bfam{e_s}_{s\in S}$ of partial  isometries such that $\sum_{s\in S}e_se_s^*=\id_E$ (strong convergence in $\sB^a(E)\subset\sB(H)$).

\item
Boundedness of amplifications of bounded right linear maps between (pre-)Hilbert \nbd{\cB}mod\-ules under tensoring from the right (Rieffel \cite[Theorem 1.4]{Rie69} or Skeide \cite[Remark 3.3]{Ske05c})%
\footnote{
The result \cite[Theorem 1.4]{Rie69} is for tensor products of Banach modules over \nbd{C^*}algebras. \cite[Lemma 3.9]{Ske00b} is for (pre-)Hilbert modules, but (as sketched in \cite[Remark 3.3]{Ske05c}; see also the proof of \cite[Proposition 1.9]{BMSS12}) the proof in \cite{Ske00b} works correctly only for pre-Hilbert modules over a von Neumann algebra; in the case of \nbd{C^*}algebras, there is a gap in the proof.  See the appendix for a correct and new proof of this important result.
}%
. In our concrete case: If a von Neumann module $E$ over $\cB\subset\sB(G)$ is represented as $E\subset\sB(G,H)$, then a bounded right linear map $\Phi\colon E\rightarrow\cB$ is represented as $\Phi\in\sB(H,G)$ via $\Phi(x)=\Phi x$ (operator multiplication); see \cite[Remark 3.3]{Ske05c}.%
\footnote{
See also the proof of \cite[Proposition 1.9]{BMSS12}.
}
Namely: If $E\subset\sB(G,H)$, then (see Footnote \ref{embFN}) $H=E\odot G$; and the general result asserts that the map $\Phi$\,``$=$''\,$\Phi\odot\id_G\colon H=E\odot G\rightarrow\cB\odot G=G$ is bounded.
\end{itemize}
With these ingredients, it is more or less plain to see that $y:=\sum_{s\in S}e_s\Phi(e_s)^*$ (strong convergence!) is an element of $E$ such that $\AB{y,x}=y^*x=\Phi(x)$.

\brem
As for these three ingredients and their interplay, polar decomposition follows (as explained in \cite{Ske00b}) directly from polar decomposition in a von Neumann algebra; it is not necessary to redo its proof for Hilbert-von Neumann modules. It is applied to make sure that the elements $e_s$ that occur in the quasi-orthonormal system (that is, $\AB{e_{s'},e_s}=\AB{e_s,e_s}\delta_{s',s}$, but not necessarily complete) can be chosen to be partial isometries. Existence of a complete such system follows by a standard application of \it{Zorn's lemma} from the observation that $\sum_{s\in S}e_se_s^*$ still converges strongly to a projection \cite[Proposition 4.9]{Ske00b}%
\footnote{
See also the proof of \cite[Proposition 1.7]{BMSS12}.
}%
).

We have already explained before this remark how boundedness of $\Phi\odot\id_G$, together with the other two ingredients, allows to conclude the proof of self-duality. We wish, however, to emphasize that while the first two ingredients are specific to the proof presented in \cite{Ske00b}%
\footnote{
See also the proof of \cite[Proposition 1.9]{BMSS12}.
}%
, it seems that establishing the third, boundedness of $\Phi\odot\id_G$, plays a decisive role in \bf{any} proof of self-duality. The proof of \cite[Proposition 6.10]{Rie74a} is based on identifying $E$ as an intertwiner space, which is easily shown to be self-dual.%
\footnote{
The proof of self-duality in \cite{Ske05c} has, in fact, a huge overlap with that of \cite[Proposition 6.10]{Rie74a}. The discussion of intertwiner spaces adds only a little to what was known already to \cite{Rie74a}. Only the proof of boundedness of $\Phi\odot\id_G$ (see also the proof of \cite[Proposition 1.9]{BMSS12}) seems to be really new -- and quite a bit simpler than that of the statement about Banach modules in \cite{Rie69} alluded to in \cite{Rie74a}.\vspace{1ex}
}
(We come back to this in Section \ref{tpSEC}.) Since the statement for (pre-)Hilbert modules in \cite[Lemma 3.9]{Ske00b} is as well-known as it is difficult to actually find a (simple and independent) proof of it%
\footnote{
It appears that many books or papers with introductions to the topic of Hilbert modules, either look only at amplifications of adjointable operators (where boundedness is an easy consequences of automatic contractivity of homomorphisms between \nbd{C^*}algebras) or look at not necessarily adjointable operators but without considering their amplifications. It is likely that Rieffel's original result in \cite{Rie69} (going beyond our scope) has been reproved in the (more general!) context of operator modules over (not necessarily self-adjoint) operator algebras.
}%
, we present a (new, we think) easy proof of this (general and interesting in its own) result in the appendix.
\erem

\bob
Here are some simple standard consequences -- \it{standard} because they are (when choosing the right proof, of course) proved exactly as for Hilbert spaces. Here, $\cB$ is a \nbd{C^*}algebra.
\begin{enumerate}
\item
Of course, a self-dual pre-Hilbert \nbd{\cB}module $E$ is necessarily Hilbert. (Corollary \ref{vNiffcor} tells that if $\cB\subset\sB(G)$ is a von Neumann algebra, then $E$ is necessarily strongly closed in the associated $\sB(G,H)$.)

\item
If $E$ is self-dual, then every bounded right linear map $E\rightarrow F$ has a (bounded right linear) adjoint.

\item
Let $F$ be a (pre-)Hilbert submodule of the pre-Hilbert \nbd{\cB}module $E$. If the submodule $F$ is self-dual then there is projection $p\in\sB^a(E)$ such that $pE=F$.%
\footnote{
See also \cite[Proposition 1.7]{BMSS12}. Of course, for strongly closed submodules of von Neumann modules, this can be proved by constructing the projection with the help of a quasi-orthonormal basis for the submodule (which is a von Neumann module in its own right). But for the simple and elegant (standard) proof from self-duality of the submodule, $\cB$ may even be only a pre-\nbd{C^*}algebra.
}

\item
And, less directly (rather a consequence of the preceding one and the fact that the strong closure of a pre-Hilbert submodule of a von Neumann module is a von Neumann module, hence, self-dual), if $A$ is a subset of a von Neumann module $E$, then the double-complement $A^{\perp\perp}$ is the von Neumann submodule of $E$ generated by $A$ (that is, $\cls^sA\cB=A^{\perp\perp}$). In particular, if $F$ is a pre-Hilbert submodule, then $F^{\perp\perp}=\ol{F}^s\subset E\subset\sB(G,H)$.

\end{enumerate}
\eob

\section{Von Neumann correspondences, Connes correspondences,\\and their tensor products} \label{tpSEC}

We do not know when exactly a full definition of the internal tensor product of two (\nbd{C^*})cor\-re\-spond\-ences%
\footnote{
Recall that a \hl{correspondence} from $\cA$ to $\cB$ (or \hl{\nbd{\cA}\nbd{\cB}correspondence}) is a Hilbert \nbd{\cB}module $E$ with a nondegenerate(\bf{!}) left action of $\cA$ that induces a homomorphism $\cA\rightarrow\sB^a(E)$.
}
occurred for the first time; but, in special cases it is present already in Rieffel's paper \cite{Rie74}. The (\hl{internal}) \hl{tensor product} (or \hl{tensor product over $\cB$}) of an \nbd{\cA}\nbd{\cB}correspondence $E$ and a \nbd{\cB}\nbd{\cC}correspondence $F$ it the unique \nbd{\cA}\nbd{\cC}correspondence $E\odot F$ generated by elements $x\odot y$ ($x\in E$, $y\in F$) subject to the relations
\baln{
\AB{x\odot y,x'\odot y'}
&
~=~
\AB{y,\AB{x,x'}y'},
&
a(x\odot y)
&
~=~
(ax)\odot y.
}\ealn
A special case is the tensor product $E\odot G$ of a Hilbert \nbd{\cB}module (that is, a \nbd{\C}\nbd{\cB}correspondence) $E$ and a representation space $G$ of $\cB$ (that is, a \nbd{\cB}\nbd{\C}correspondence), we have already met (for faithful representations of $\cB$ only, but faithfulness is, of course, not necessary).%
\footnote{ \label{tpFN}
By simple \nbd{C^*}algebraic arguments it follows that the sesquilinear map on the algebraic tensor product $E\otimes F$ defined as in the first relation above, is positive, hence, defines a semiinner product on $E\otimes F$. (See the appendix.) Doing Hausdorff completion we obtain the Hilbert \nbd{\cC}module $E\odot F$, and since the left action of $\cA$ induced by the second relation is by adjointable elements, it survives the quotient.

However, already in \cite[Section 4.2]{Ske01} the following alternative has been observed: It is easy, using cyclic decomposition of the representation of $\cB$ on $G$ (see also the proofs \cite[Lemma 3.9]{Ske00b} and \cite[Proposition 1.9]{BMSS12}), to see positivity of the inner product of $E\odot G$. Using this, (assuming $\cC\subset\sB(L)$) one gets that $x\odot y$ may be represented as operator $L_xL_y\in\sB(L,E\odot(F\odot L))$ where $L_y=y\odot\id_L\in\sB(L,F\odot L)$ and $L_x=x\odot\id_{F\odot L}\in\sB(F\odot L,E\odot(F\odot L))$; see, in particular, \cite[Equation (4.2.2)]{Ske01}. (See also the discussion leading to \cite[Proposition 3.4]{BMSS12}.)
}

A (\hl{concrete}) \hl{von Neumann correspondence}~ (\cite[Definition 4.2]{Ske06b}) ~\cite[Definition 3.3.1]{Ske01} is a correspondence between von Neumann algebras which is a (concrete) von Neumann module such that the left action is \it{normal}. (There are various equivalent ways, to say what it means that the left action is normal; see, for instance, \cite[Lemma  3.3.2]{Ske01}. Just recall that a (necessarily positive\bf{!}) map is \hl{normal} if it is order-continuous and that, for positive maps, this is equivalent to \nbd{\sigma}weak continuity.)

\brem
There is an -- obvious -- analogue of Theorem \ref{=thm} for Hilbert-von Neumann bimodules \cite{BMSS12} and concrete von Neumann correspondences \cite{Ske06b}, which we omit stating. And, obviously, every concrete von Neumann correspondence is a von Neumann correspondence, while every von Neumann correspondence can be turned into a concrete one be the procedure mentioned in Footnote \ref{embFN}.
\erem

Let $\cA\subset\sB(K)$, $\cB\subset\sB(G)$, $\cC\subset\sB(L)$ be von Neumann algebras and let $E$ and $F$ be a von Neumann \nbd{\cA}\nbd{\cB}correspondence and a von Neumann \nbd{\cB}\nbd{\cC}correspondence, respectively. Then:
\begin{itemize}
\item
Either construct the internal tensor product $E\odot F$, embed it into $\sB(L,(E\odot F)\odot L)$ (see Footnote \ref{embFN}) to obtain the strong closure $E\sodots F:=\cls^sE\odot F\subset\sB(L,E\odot F\odot L)$.

\item
Or, construct for $x\in E$ and $y\in F$ the operators $L_y=y\odot\id_L\in\sB(L,F\odot L)$ and $L_x=x\odot\id_{F\odot L}\in\sB(F\odot L,E\odot(F\odot L))$ (see Footnote \ref{tpFN}) to obtain the strong closure $E\sodots F:=\cls^sL_EL_F\subset\sB(L,E\odot F\odot L)$.
\end{itemize}
Since in the tensor category of correspondences the spaces $(E\odot F)\odot L$ and $E\odot(F\odot L)$ are identified and denoted as $E\odot F\odot L$, the two possibilities coincide. So the two ways lead to the same von Neumann \nbd{\cC}module $E\sodots F$. It is a remarkable feature of von Neumann modules obtained from strong closures of (pre-)Hilbert modules over a von Neumann algebra, that the bounded right linear operators on the original (pre-)Hilbert module extend as (adjointable) operators on the strong closure. (Taking into account the statement proved in the appendix, the bounded right linear operators are represented as bounded operators between the associated Hilbert spaces, and the composition by operator multiplication extends, of course, to the strong closure.) It is routine to show that the left action of $\cA$ on $E\odot F$ extended to $E\sodots F$ is normal. So, $E\sodots F$ is a von Neumann \nbd{\cA}\nbd{\cC}correspondence, the (\hl{internal}) \hl{tensor product} in the category of von Neumann correspondences.

\brem
In \cite[Section 3]{BMSS12}, the authors call the (internal) tensor product of Hilbert-von Neumann bimodules (that is, of concrete von Neumann correspondences) the `\it{Connes fusion}'. Reasoning: There are certain bimodules over von Neumann algebras (also known under the name \it{Connes correspondences}) and a sort of tensor product for them (sometimes going under the name of `\it{Connes fusion}')%
\footnote{
We have heard rumors that the creator of that construction is not too fond of that terminology, so we shall not use it.
}%
, so it would be natural to use the same terminology also for other bimodules \it{of von Neumann type}.
\begin{itemize}
\item
We find this a rather feeble motivation for changing the name of something, the tensor product over an algebra, that existed quite a bit before Connes correspondences.

\item
The motivation does not right hit the point. In fact, as mentioned in \cite[Remark 4.3]{Ske06b}, Connes correspondences may be viewed as (a special case of) concrete von Neumann correspondences and, in this translation, the tensor product of the latter \bf{is} Connes' product of the former. These facts are entirely missing in \cite{BMSS12}.%
\footnote{
We did not quite understand, if no relation between Hilbert-von Neumann bimodules and Connes correspondences is made in \cite{BMSS12}, or if the authors think the relation would be the occurrence of two representations in the definition of Hilbert-von Neumann bimodule. The latter would not be applicable, because for a Hilbert-von Neumann bimodule $E\subset\sB(G,H)$ there is the representation $\pi$ as for Hilbert-von Neumann modules (see the beginning of Section \ref{compSEC}) acting on $G$ plus an analogue representation acting on $H$, while a Connes correspondence (see the basic fact 5 below) is a single Hilbert space $H$ with a representation and an anti-representation both acting on $H$.\vspace{1ex}
}%

We shall explain and prove this in the sequel.
\end{itemize}
\erem

\noindent
A concrete von Neumann correspondence $E\subset\sB(G,H)$ from a von Neumann algebra $\cA\subset\sB(K)$ to a von Neumann algebra $\cB\subset\sB(G)$ comes along with a normal and unital representation
\beqn{
\rho
\colon
\cA
~\longrightarrow~
\sB^a(E)
~\subset~
\sB(H)
}\eeqn
(where, with $H=E\odot G$, we identify $\sB^a(E)=\sB^a(E)\odot\id_G\subset\sB(H)$; see Footnote \ref{embFN}), which we call the \hl{Stinespring representation}%
\footnote{
Already in \cite[Section 7]{Ske00b} (see also \cite[Example 2.16 and Observation 2.17]{BhSk00}) it has been explained that when $E$ is the \it{GNS-correspondence} (Paschke \cite{Pas73}) of a CP-map into a concrete \nbd{C^*}algebra $\cB\subset\sB(G)$, then $\rho$ is, indeed, the Stinespring representation. See also \cite[Section 2]{BMSS12}.\vspace{1ex}
}
of $\cA$ on $H$ associated with $E$. The Stinespring representation comes from amplification of operators in $\sB^a(E)$ acting on the left tensor factor $E$ in the tensor product $H=E\odot G$. However, there is also the possibility to amplify operators acting on the right factor $G$. But, in order that an operator $b'$ can by amplified to an operator $\id_E\odot b'$ on $H$, it has to be left \nbd{\cB}linear, that is exactly, $b'\in\cB'\subset\sB(G)$. We refer to this normal and unital representation $\rho'\colon\cB'\rightarrow\id_E\odot\cB'\subset\sB(H)$ as the \hl{commutant lifting}%
\footnote{
Arveson's representation discussed in the section ``\it{Lifting Commutants}'' in \cite{Arv69}, is a special case of our \it{commutant liftings}.
}
of $\cB'$ on $H$ associated with $E$.

We collect some basic facts:
\begin{enumerate}
\item \label{C1}
$E=C_{\cB'}(\sB(G,H)):=\bCB{x\in\sB(G,H)\colon\rho'(b')x=xb'~(b'\in\cB')}$.

\item\label{C2}
Conversely, if $\rho'$ is any normal unital representation of $\cB'$ on a Hilbert space $H$, then $E\subset\sB(G,H)$ defined as $C_{\cB'}(\sB(G,H))$ is a concrete von Neumann \nbd{\cB}module and its commutant lifting is the $\rho'$ we started with. \cite[Theorem 3.10]{Ske06b}: The category of concrete von Neumann \nbd{\cB}modules and the category of normal unital representations of $\cB'$ are the same (not just equivalent, but really ``the same'', in the sense that objects and morphisms correspond one-to-one).

\item\label{C3}
$\rho'(\cB')'=\sB^a(E)$. Since $\rho(\cA)\subset\sB^a(E)$, we see that, in particular, $\bSB{\rho(\cA),\rho(\cB')}=\zero$.

\item\label{C4}
Conversely, if $\rho$ and $\rho'$ are any commuting pair of normal unital representations of $\cA$ and $\cB'$, respectively, then $E:=C_{\cB'}(\sB(G,H))\subset\sB(G,H)$ is a concrete von Neumann \nbd{\cB}module inheriting a left action of $\cA$ via $\rho$. Also here, we have a true one-to-one correspondence between the categories of such commuting pairs and concrete von Neumann \nbd{\cA}\nbd{\cB}correspondences.

But, we have even more. Exchanging the roles of $\rho$ and $\rho'$, we get a concrete von Neumann \nbd{\cB'}\nbd{\cA'}correspondence as $E'=C_\cA(\sB(K,H)):=\bCB{x'\in\sB(K,H)\colon\rho(a)x'=x'a~(a\in\cA)}$ with left action of $\cB'$ via $\rho'$. Also here we have \cite[Theorem 4.4]{Ske06b}: This \hl{commutant} is a bijective functor between the category of concrete von Neumann \nbd{\cA}\nbd{\cB}correspondences and the category of concrete von Neumann \nbd{\cB'}\nbd{\cA'}correspondences.

\item\label{C5}
In particular, if $\cB$ is in \it{standard representation}, then $\cB'\cong\cB^\circ$, the opposite von Neumann algebra of $\cB$, via \it{Tomita conjugation}. That is, on $H$ we have a commuting pair of representations $ \rho$ of $\cA$ and $\rho^\circ$ of $\cB^\circ$. In other words, $H$ is an \nbd{\cA}\nbd{\cB}bimodule, and as such a \hl{Connes correspondence} from $\cA$ to $\cB$; see Connes \cite{Con80p}. Conversely, interpreting the representation $\rho^\circ$ of $\cB^\circ$ as representation $\rho'$ of $\cB'$, every Connes correspondence gives rise to a concrete von Neumann correspondence.

But, in general, for us $\cB$ need not be in standard representation: Every Connes correspondence is a concrete von Neumann correspondence, but not every concrete von Neumann correspondence is a Connes correspondence.
\end{enumerate}

\brem
The proofs of \eqref{C1} and \eqref{C3} rely on the simple computation
\beqn{
\SMatrix{\cB&E^*\\E&\sB^a(E)}''
~=~
\left\{\SMatrix{b'&\\&\rho'(b')}\colon b'\in\cB'\right\}'
~=~
\SMatrix{\cB&C_{\cB'}(\sB(E\odot G,G))\\C_{\cB'}(\sB(G,E\odot G))&\rho'(\cB')'}
}\eeqn
(which holds for every pre-Hilbert module over $\cB\subset\sB(G)$, when embedded into $\sB(G,H)$ with $H=E\odot G$); see \cite[Section 2]{Ske05c}. The proof of \eqref{C2}, hence, \eqref{C4} (only nondegeneracy $\cls EG=H$ is an issue), relies, furthermore, on Muhly and Solel \cite[Lemma 2.10]{MuSo02}. Both have been first observed in \cite[Section 2]{Ske03c}. The statement in \eqref{C5} is not more than a simple observation that has been made in \cite[Remark 4.3]{Ske06b}.
\erem

\bthm
The internal tensor product of two Connes correspondences (viewed as von Neumann correspondences) is (when viewed as Connes correspondence) Connes' product of the two Connes correspondences.
\ethm

This has been mentioned after Skeide \cite[Remark 6.3]{Ske08a} without indicating a proof. The argument given in Shalit and Skeide \cite[Footnote y]{ShaSk10p} is very sketchy. Therefore:

\proof
Let $\cA$, $\cB$, $\cC$ be \nbd{W^*}algebras, let $H$ be a Connes \nbd{\cA}\nbd{\cB}correspondence, and let $K$ be a Connes \nbd{\cB}\nbd{\cC}correspondence. As far as their Connes product $H\otimes_\psi K$ (relative to a faithful normal semi-finite weight $\psi$ on $\cB$) is concerned, we use the presentation in Takesaki \cite[Section IX.3]{Tak03a}.

For the time being, we do not fix a representation of $\cA$ nor of $\cC$ turning them into von Neumann algebras. (This is only necessary in the very end, when we wish to consider $H$ and $K$ as von Neumann correspondences which are defined only over von Neumann algebras.) But we assume that $\cB\subset\sB(G)$ where $G=L^2(\cB,\psi)$, so that $\cB$ is in standard representation and $\cB'=\cB^\circ$ via \it{Tomita conjugation}.

It follows that $H$ carries a normal  unital representation of $\cB'$ (commuting with the action of $\cA$) so that we may define the concrete von Neumann \nbd{\cB}module $E:=C_{\cB'}(\sB(G,H))$ (a von Neumann \nbd{\cA}\nbd{\cB}correspondence as soon as we represent $\cA$ as a von Neumann algebra). Note that this $E$ is exactly what has been denoted as $\sL(L^2(\sN)_\sN,\eH_\sN)$ preceding \cite[Lemma IX.3.3]{Tak03a} for $\sN=\cB$, $L^2(\sN)=G$, and $\eH=H$.

There is a dense (\cite[Lemma IX.3.3(iii)]{Tak03a}) subspace $\sD$ (denoted $\eD(\eH,\psi)$ and defined in \cite[Equation IX.3.(10)]{Tak03a}) and an injective (\cite[Lemma IX.3.3(iv)]{Tak03a}) map $\eL\colon\sD\rightarrow E$ (denoted $L_\psi$ and defined in \cite[Equation IX.3.(11)]{Tak03a}). From the definition \cite[Equation IX.3.(11)]{Tak03a}, it follows that $\eL(H)G$ is dense in $\sD$, hence, in $H$. Therefore, $\eL(\sD)$ is a strongly dense subset of $E$.

Putting in \cite[Proposition IX.3.15]{Tak03a} $\eK=K$, just looking at \cite[Equation IX.3.(23)]{Tak03a}, tells us that the (semi)inner product of elements $d_i\otimes k_i$ in $\sD\,\ul{\otimes}\,K$ defined there is
\beqn{
B(d_1\otimes k_1,d_2\otimes k_2)
~=~
\AB{\eL(d_1)\odot k_1,\eL(d_2)\odot k_2},
}\eeqn
the same as the inner product of the elements $\eL(d_i)\odot k_i$ in $E\odot K$. Now, if also $\cC\subset\sB(L)$ is in standard representation (so that $\cC'=\cC^\circ$), then $F:=C_{\cC'}(\sB(L,K))$ is a von Neumann \nbd{\cB}\nbd{\cC}correspondence such that $K=F\odot L$, including the representations of $\cB$ and and $\cC'$. It is routine to show that $\eH\otimes_\psi\eK$ in \cite[Definition IX.3.16]{Tak03a} is $E\odot F\odot L$ including the representations of $\cA$ and $\cC^\circ=\cC'$ that identify $E\sodots F$ as the von Neumann \nbd{\cA}\nbd{\cC}correspondence coming from the Connes \nbd{\cA}\nbd{\cC}correspondence $\eH\otimes_\psi\eK$.\qed

\brem
The commuting pair picture of von Neumann correspondences is, therefore, a generalization of Connes correspondences and the tensor product of the former includes the Connes product of the latter. It coincides if we only look at von Neumann algebras in standard representation.  But, unlike the construction of  Connes' product, the construction of the tensor product of von Neumann correspondences does not depend on more advanced tools like \it{Tomita-Takesaki theory} and \it{standard representation} and the formulae do not depend manifestly on the choice of a faithful semi-finite normal weight.

This has been mentioned after \cite[Remark 6.3]{Ske08a}; it has nothing to do with the term `\it{Connes fusion}' in \cite{BMSS12}. See also Thom \cite{Tho11}, where the discussion is, however, limited to \it{bi-finite} Connes correspondences (that is, the module dimensions of $H$ as left and right module are finite). We do not have any limitation.

We also would like to mention that there is a close connection with Sauvageot's \it{tensor product relative to an algebra} \cite{Sau80,Sau83}.
\erem

% \newpage

\section*{Appendix} \stepcounter{section} \renewcommand{\thesection}{A} \setcounter{emp}{0}

The scope of this appendix is to present a new and `painless' proof of the following result.

\bthm
Let $E$ be a (pre-)Hilbert module over a (pre-)\nbd{C^*}algebra $\cB$ and let $F$ be a (pre-) \nbd{C^*}correspondence from $\cB$ to another (pre-)\nbd{C^*}algebra $\cC$. Then for every $a\in\sB^r(E)$ we have
\beqn{
\norm{a\odot\id_F}
~\le~
\norm{a},
}\eeqn
so $a\odot\id_F\in\sB^r(E\odot F)$.
\ethm

\bob
{~}
\begin{itemize}
\item
We have not yet specified what we mean by $E\odot F$ in case of `pre-'. It is simply $\ls\bCB{x\odot y\colon x\in E,y\in F}\subset\ol{E}\odot\ol{F}$.%
\footnote{ \label{algtpFN}
Note that the algebraic tensor product of a right \nbd{\cB}module $E$ and a left \nbd{\cB}module $F$ \hl{over $\cB$} is $E\,\ul{\odot}_\cB F:=(E\,\ul{\otimes}\,F)/\fls\CB{xb\otimes y-x\otimes by}$. In the setting of the theorem, the semiinner product defined on $E\,\ul{\otimes}\,F$ as usual always survives the quotient but still need not be inner on $E\,\ul{\odot}_\cB F$. (It is inner if $\cB$ and $E$ are complete; see \cite[Proposition 4.2.22]{Ske01}.) Therefore, it is pointless to divide first out the relations $xb\otimes y-x\otimes by$ and then a still possible kernel of the semiinner product on $E\,\ul{\odot}_\cB F$, but we divide out immediately the kernel of the semiinner product on $E\,\ul{\otimes}\,F$. We emphasize, however, that the relation $xb\odot y=x\odot by$ is satisfied at any stage.
}
So, since $a$ is bounded we simply may complete the spaces $E$, $F$, $\cB$, $\cC$, proving the theorem only without `pre-'.

\item
Note that we need not have equality. Effectively, $E\odot F$, hence $a\odot\id_F$, may very well be zero independently of $\norm{a}$.

\item
Of course, the theorem remains true for $a\in\sB^r(E_1,E_2)$; simply embed $\sB^r(E_1,E_2)$ into $\sB^r\rtMatrix{E_1\\E_2}$.
\end{itemize}
\eob

\brem
For Banach modules over \nbd{C^*}algebras, the theorem is \cite[Theorem 1.4]{Rie69}. There should be versions for \it{operator modules} over (non-selfadjoint) \it{operator algebras} (including similar statements for CB-norms). The theorem is easy, if $a$ is required adjointable. ($a\mapsto a\odot\id_F$ is a homomorphism between \nbd{C^*}algebras and, therefore, contractive.)
\erem

The proof of the theorem is a typical example for the following strategy to prove statements about (internal) tensor products.
\bob
{~}
\begin{enumerate}
\item
Frequently, a certain statement is easy (or easier) to prove on simple tensors $x\odot y$ instead of general elements $\sum_{i=1}^nx_i\odot y_i$ of the algebraic tensor product.

\item
Use the equality of $E\odot F=E_n\odot F^n$ via $\sum_{i=1}^nx_i\odot y_i=X_n\odot Y^n$ (with $X_n=(x_1,\ldots,x_n)\in E_n$ and $Y^n=\rtMatrix{y_1\\\vdots\\y_n}\in F^n$, where $E_n$ is the row-module with inner product $\AB{X_n,X'_n}=\bfam{\AB{x_i,x'_j}}\in M_n(\cB)$ and the obvious action of $M_n(\cB)$ on $F^n$) to reduce the problem to simple tensors.
\end{enumerate}
Among the various applications, there is that the inner product on $E\odot F$ is actually positive ($\BAB{\sum_{i=1}^nx_i\odot y_i,\sum_{i=1}^nx_i\odot y_i}=\AB{X_n\odot Y^n,X_n\odot Y^n}=\AB{Y^n,\AB{X_n,X_n}Y^n}\ge0$, having established before that the inner product of $E_n$ is positive; see Skeide \cite[Section 5]{Ske11}). Another application is the proof in Skeide and Sumesh \cite[Corollary]{SkSu14} of Blecher's result \cite[Theorem 4.3]{Ble97} that $E\odot F$ is (completely isometrically) isomorphic to the Haagerup tensor product.
\eob

\proof[Proof of the theorem.]
Recall that we may restrict to the case without `pre-'. The statement of the theorem actually means that we show the estimate for elements of the form $\sum_{i=1}^nx_i\odot y_i$ and, then, extend extend $a\odot\id_F$, well defined on that dense subspace, continuously to all of $E\odot F$. (Recall from Footnote \ref{algtpFN} that this subspace is $E\,\ul{\odot}_\cB F$.)

The proof decomposes into four parts.

\item
For all $a\in\sB^r(E)$ we have $\AB{ax,ax}\le\norm{a}^2\AB{x,x}$ in $\cB$. (This is the easier direction of \cite[Theorem 2.8]{Pas73}, which asserts that the inequality actually characterizes among all linear maps on $E$ those that are in $\sB^r(E)$. See there for the quick proof of the part we need here; the other direction is much harder.)

\item
For each $a\in\sB^r(E)$ we have $\sB^a(E)\ni a\ge0$ if (and only if) $\AB{x,ax}\ge0$ for all $x\in E$. (In this form, this is Lance \cite[Lemma 4.1]{Lan95}, while \cite[Proposition 6.1]{Pas73} is the special case we would actually need in the next part. Paschke actually shows in his proof Lance' more general statement; but here we prefer Lance' proof.)

\item
The inequality in Part (i) survives `amplification' to $E_n$. More precisely, for each $X_n=(x_1,\ldots,x_n)\in E_n$, we have $\Bfam{\AB{ax_i,ax_j}}\le\norm{a}^2\Bfam{\AB{x_i,x_j}}$ in $M_n(\cB)$ or, equivalently,
\beqn{
M_n(\cB)
~\ni~
\Bfam{\norm{a}^2\AB{x_i,x_j}-\AB{ax_i,ax_j}}
~\ge~0.
}\eeqn
Since $M_n(\cB)$ acts in the obvious way on the Hilbert \nbd{\cB}module $\cB^n$, we may apply Part (ii) to see  that this element of $M_n(\cB)$ is positive. Indeed, for $B^n=\rtMatrix{n_1\\\vdots\\b_n}\in\cB^n$ we find
\baln{
\AB{B^n,\Bfam{\norm{a}^2\AB{x_i,x_j}-\AB{ax_i,ax_j}}B^n}
&
~=~
\sum_{i,j}b_i\bfam{\norm{a}^2\AB{x_i,x_j}-\AB{ax_i,ax_j}}b_j
\\
&
~=~
\textstyle\norm{a}^2\BAB{\sum_ix_ib_i,\sum_ix_ib_i}-\BAB{a\sum_ix_ib_i,a\sum_ix_ib_i}.
}\ealn
By Part (i), this is positive.

\item
By Part (iii) (applied in the step where the $\le$ occurs), for $Y^n=\rtMatrix{y_1\\\vdots\\y_n}\in F^n$ we find
\bmun{
\Bnorm{(a\odot\id_F)\Bfam{\sum_ix_i\odot y_i}}^2
~=~
\Bnorm{\Bfam{\sum_iax_i\odot y_i}}^2
~=~
\norm{aX_n\odot Y_n}^2
~=~
\bnorm{\AB{Y^n,\AB{aX_n,aX_n}Y^n}}
\\
~\le~
\norm{a}^2\bnorm{\AB{Y^n,\AB{X_n,X_n}Y^n}}
~=~
\norm{a}^2\Bnorm{\sum_ix_i\odot y_i}^2,
}\emun
what was to be shown.\qed

\lf\lf\noindent
\bf{Acknowledgments.} Discussions with several people helped a lot to improve this paper. In particular, I would like to thank Malte Gerhold and Orr Shalit.

% \newpage

\setlength{\baselineskip}{2.5ex}

% \vspace{-2ex}
% \bibliography{mybib}
\newcommand{\Swap}[2]{#2#1}\newcommand{\Sort}[1]{}
\providecommand{\bysame}{\leavevmode\hbox to3em{\hrulefill}\thinspace}
\providecommand{\MR}{\relax\ifhmode\unskip\space\fi MR }
\providecommand{\MRhref}[2]{%
  \href{http://www.ams.org/mathscinet-getitem?mr=#1}{#2}
}
\providecommand{\href}[2]{#2}

\lf\noindent
Michael Skeide: \it{Dipartimento di Economia, Università degli Studi del Molise, Via de Sanctis, 86100 Campobasso, Italy},
E-mail: \href{mailto:skeide@unimol.it}{\tt{skeide@unimol.it}},
Homepage: \url{http://web.unimol.it/skeide}

\end{document}